\documentclass[12pt]{amsart}
\usepackage{amssymb, enumerate}
\theoremstyle{plain}

\newtheorem{thm}{Theorem}[section]
\newtheorem{cor}[thm]{Corollary}
\newtheorem{lem}[thm]{Lemma}
\newtheorem{prop}[thm]{Proposition}

\theoremstyle{definition}
\newtheorem{defi}[thm]{Definition}
\newtheorem{conj}[thm]{Conjecture}
\newtheorem{conv}[thm]{Convention}
\newtheorem{nota}[thm]{Notation}
\newtheorem{rem}[thm]{Remark}
\newtheorem{rems}[thm]{Remarks}
\newtheorem{exa}[thm]{Example}
\newtheorem{exas}[thm]{Examples}
\newtheorem{sit}[thm]{}

\newcommand{\brem}{\begin{rem}}
\newcommand{\brems}{\begin{rems}}
\newcommand{\erem}{\end{rem}}
\newcommand{\erems}{\end{rems}}
\newcommand{\bexa}{\begin{exa}}
\newcommand{\bexas}{\begin{exas}}
\newcommand{\eexa}{\end{exa}}
\newcommand{\eexas}{\end{exas}}
\newcommand{\bdefi}{\begin{defi}}
\newcommand{\edefi}{\end{defi}}
\newcommand{\bcor}{\begin{cor}}
\newcommand{\ecor}{\end{cor}}
\newcommand{\blem}{\begin{lem}}
\newcommand{\elem}{\end{lem}}
\newcommand{\bconv}{\begin{conv}}
\newcommand{\econv}{\end{conv}}
\newcommand{\bconj}{\begin{conj}}
\newcommand{\econj}{\end{conj}}
\newcommand{\bprop}{\begin{prop}}
\newcommand{\eprop}{\end{prop}}
\newcommand{\bthm}{\begin{thm}}
\newcommand{\ethm}{\end{thm}}
\newcommand{\bnota}{\begin{nota}}
\newcommand{\enota}{\end{nota}}
\newcommand{\bsit}{\begin{sit}}
\newcommand{\esit}{\end{sit}}
\newcommand{\be}{\begin{eqnarray}}
\newcommand{\ee}{\end{eqnarray}}
\newcommand{\bproof}{\begin{proof}}
\newcommand{\eproof}{\end{proof}}

\def\ba{\begin{array}}
\def\ea{\end{array}}

\def\cO{{\mathcal O}}

\newcommand{\Spec}{\operatorname{Spec}}
\newcommand{\Pic}{\operatorname{Pic}}
\newcommand{\Cl}{\operatorname{Cl}}
\newcommand{\Frac}{\operatorname{Frac}}

\newcommand{\Proj}{\operatorname{Proj}}

\newcommand{\id}{\operatorname{id}}

\newcommand{\Aut}{{\operatorname{Aut}}}
\newcommand{\ML}{{\operatorname{ML}}}
\newcommand{\LND}{{\operatorname{LND}}}
\newcommand{\rk}{{\operatorname{rk}}}
\newcommand{\SL}{{\bf {SL}}}
\newcommand{\GL}{{\bf {GL}}}

\def\PGL{{\bf PGL}}

\renewcommand{\div}{{\operatorname{div}}}

\newcommand{\A}{{\mathbb A}}
\newcommand{\PP}{{\mathbb P}}

\newcommand{\C}{{\mathbb C}}
\newcommand{\Q}{{\mathbb Q}}
\newcommand{\Z}{{\mathbb Z}}
\newcommand{\N}{{\mathbb N}}
\newcommand{\T}{{\mathbb T}}

\newcommand{\p}{{\partial}}

\def\tdelta{{\tilde\delta}}

\title[Uniqueness of $\C^*$-actions]{On the uniqueness of $\C^*$-actions
on affine surfaces}

\author{Hubert Flenner}
\address{Fakult\"at f\"ur Mathematik,
Ruhr Universit\"at Bochum,
Geb.\ NA 2/72,
Universit\"ats\-str.\ 150,
44780 Bochum, Germany}
\email{Hubert.Flenner@ruhr-uni-bochum.de}

\author{Mikhail Zaidenberg}
\address{Universit\'e
Grenoble I, Institut Fourier, UMR 5582 CNRS-UJF, BP 74,
38402 St.\ Martin
d'H\`eres c\'edex, France}
\email{zaidenbe@ujf-grenoble.fr}

\thanks{
{\bf Acknowledgements:} This research was done during a visit of
the first author at the Institut Fourier of the University of
Grenoble. He thanks this institution for the generous support and
excellent working conditions.}

\thanks{
\mbox{\hspace{11pt}}{\it 1991 Mathematics Subject
Classification}:
14R05, 14R20, 14J50.\\
\mbox{\hspace{11pt}}{\it Key words}: $\C^*$-action, $\C_+$-action,
graded algebra, affine surface}

\begin{document}

\begin{abstract}
It is an open question whether every normal affine surface $V$
over $\C$ admits an effective action of a maximal torus
$\T=\C^{*n}$ ($n\le 2$) such that any other effective
$\C^*$-action is conjugate to a subtorus of $\T$ in $\Aut (V)$. We
prove that this holds indeed in the following cases: (a) the
Makar-Limanov invariant $\ML(V) \ne\C$ is nontrivial, (b) $V$ is
a toric surface, (c) $V=\PP^1\times\PP^1\backslash \Delta$, where
$\Delta$ is the diagonal, and (d) $V=\PP^2\backslash Q$, where
$Q$ is a nonsingular quadric. In case (a) this generalizes a
result of Bertin for smooth surfaces, whereas (b) was previously
known for the case of the affine plane (Gutwirth \cite{Gut}) and
(d) is a result of Danilov-Gizatullin \cite{DG} and Doebeli
\cite{Do}.
\end{abstract}

\maketitle

\section{Introduction}

The classification problem for reductive group actions on affine
spaces or, more generally, on affine varieties, has a long
history. By \cite{Kam, KP, KKMLR, Po} any reductive group action
on $\A^2_\C$ and $\A^3_\C$ is conjugate to a linear one. The same
holds for connected reductive groups acting on $\A^4_\C$ except
possibly for $\C^*$ and $\C^{*2}$ \cite{Pa, Po} (cf. also
\cite{BaHab}), and for tori $\T^n$ acting effectively on
$\A_\C^n$, $\A_\C^{n+1}$ \cite{BB} and on the affine toric
$n$-folds \cite{De, Gub}. According to \cite{Sc, MMP, Kn} many
finite nonabelian groups and any connected reductive nonabelian
group admit a non-linearizable action on some affine space
$\A^n_\C$. In the local case the existence of a maximal reductive
subgroup of $\Aut(V,0)$, which contains a conjugate of any other
connected reductive subgroup, was established in \cite{HM}. In
\cite{DG, Do} the same was shown to be true for the smooth affine
quadric surface in $\A^3_\C$. It is an open question whether this
holds as well for every normal affine surface $V$. In this paper
we give some partial positive results, see Corollary \ref{fc}
below.

Bertin's Theorem \cite[Corollary 2.3]{Be} asserts that, for a
smooth affine surface $V$ non-isomorphic to $\C^*\times \C^*$,
which admits a minimal compactification $V$ by a simple normal
crossing divisor $D$ with a non-linear dual graph $\Gamma_D$, any
two effective $\C^*$-actions on $V$ are conjugated in the
automorphism group $\Aut(V)$. On the other hand, by Gizatullin's
Theorem \cite[Theorems 2 and 3]{Gi} (see also \cite{BML, Be} or
\cite{Du} for the more general case of normal surfaces), if
$V\not\cong \A^1_\C\times \C^*$ then $\Gamma_D$ is linear if and
only if the Makar-Limanov invariant
$$\ML(V):=\bigcap_{\p\in\LND(A)} \ker \p$$
of $V$ is trivial that is, $\ML(V)=\C$, where $\LND(A)$ stands
for the set of all locally nilpotent derivations of the
coordinate ring $A=H^0(V,\cO_V)$ of $V$. The latter holds if and
only if $V$ admits two non-equivalent effective $\C_+$-actions
i.e., two $\C_+$-actions with different general orbits (see e.g.,
\cite{FZ2}).

We present here an alternative proof of Bertin's Theorem, valid
more generally for normal affine surfaces. Our proof is not based
on the properties of completions and so is independent of
Gizatullin's Theorem. In Theorem \ref{0.6} below we show that, as
soon as $\ML(V)\neq\C$ and $V\not\cong \C^*\times\C^*$ or
$\A_\C^1\times\C^* $, any two effective $\C^*$-actions on $V$ are
conjugated via an element of a $\C_+$-subgroup of the
automorphism group $\Aut (V)$.

For a surface $V$ with an effective action of the 2-torus $\T$
we prove in Theorem \ref{4.2} that any effective $\C^*$-action on
$V$ is conjugated in $\Aut (V)$ to the action of a subtorus of
$\T$. In the case of the affine plane $V=\A_\C^2$ this gives
another proof of the classical Gutwirth Theorem \cite{Gut} saying
that the linearization conjecture for $\C^*$-actions on $\A_\C^n$
holds in dimension 2.

In Section 5 we deduce similar results for the surfaces
$\PP^1\times\PP^1\backslash\Delta$ and $\PP^2\backslash Q$, where
$\Delta$ is the diagonal in $\PP^1\times\PP^1$ and $Q$ is a
nonsingular quadric in $\PP^2$. We show that any two effective
$\C^*$-actions on one of these surfaces are conjugate in the
automorphism group. According to a result of Gizatullin and Popov
(see \cite[Theorem 4.12]{FZ2}) these are the only normal affine
surfaces that admit a nontrivial $\SL_2$-action except for the
affine plane and the affine Veronese cones over the rational
normal curves. We note that by \cite[Proposition 4.14]{FZ2} any
two $\SL_2$-actions on a normal affine surface are conjugate in
the automorphism group.

Our interest in such kind of results is related with our studies
\cite{FZ1, FZ2} on the Dolgachev-Pinkham-Demazure (or DPD, for
short) presentation of a normal affine surface $V$ endowed with a
$\C^*$-action. We show in Corollary \ref{4.1} that for surfaces
with a non-trivial Makar-Limanov invariant, except in the case of
the surfaces $\C^*\times\C^*$ and $\A^1_C\times \C^*$, this
DPD-presentation is uniquely determined up to a natural
equivalence. For affine toric surfaces we describe in Section 4
the possible ambiguities in the choice of a DPD-presentation. We
also deduce the uniqueness of the DPD-presentation for the
surfaces $\PP^1\times\PP^1\backslash\Delta$ and $\PP^2\backslash
Q$ as above, which are non-toric and have a trivial Makar-Limanov
invariant.

In Corollary \ref{fc} we deduce that all maximal connected
reductive subgroups of the automorphism group $\Aut (V)$ are
conjugate, and any connected reductive subgroup is contained in a
maximal one, besides the remaining open case when the surface is
non-toric and has a trivial Makar-Limanov invariant. In the
forthcoming paper \cite{FKZ} we will solve this remaining case,
up to one exception, by  showing that the automorphism group $\Aut
(V)$ contains a unique class of conjugated $\C^*$-subgroups, which
also implies the uniqueness of a DPD-presentation up to a natural
equivalence.

\section{Preliminaries }

Let $k$ be an algebraically closed field of characteristic $0$. We
recall the following definition.

\bdefi\label{0.1} Let $A$ be a $k$-algebra, not necessarily
associative or commutative. A derivation $\delta:A\to A$ is
called {\em locally bounded} if every finite subset of $A$ is
contained in some $\delta$-invariant linear subspace $V\subseteq
A$ of finite dimension over $k$.

For instance, if $\delta:A\to A$ is a semisimple derivation on
$A$ i.e., $A$ has a $k$-basis consisting of eigenvectors of
$\delta$ then $\delta$ is locally bounded. \edefi

Clearly if $\delta$ is locally bounded then $A$ is the union of
its finite dimensional $\delta$-invariant subspaces $V$. In
particular, the restriction of $\delta$ to each such subspace $V$
admits a unique Jordan-Chevalley decomposition
$$
\delta=\delta_s+\delta_n,
$$
where $\delta_s$ and $\delta_n$ are the semisimple and the
nilpotent parts of $\delta$, respectively. The restriction to a
smaller $\delta$-invariant subspace respects this decomposition,
hence $\delta_s$ and $\delta_n$ are well-defined $k$-linear maps
on $A$ that are again locally bounded.  We will see in Lemma
\ref{0.2} below that these maps are as well derivations on $A$.

For $\alpha\in k$ we consider the linear subspace
$$
A_\alpha:=\bigcup_{i\ge 0}\ker (\delta-\alpha\id)^i\,.
$$
   From standard linear algebra we know that $$V=\bigoplus_{\alpha}
(V\cap A_\alpha)\,,$$ whenever $V$ is a finite dimensional
$\delta$-invariant subspace of $A$. Thus $A$ is a graded vector
space
$$
A=\bigoplus_{\alpha\in k}A_\alpha\,.
$$ Clearly, $\delta_s$ and $\delta_n$ leave invariant every
subspace $A_\alpha$. Moreover $\delta_s$ acts on $A_\alpha$ via
multiplication by $\alpha$, whereas $\delta_n|A_\alpha$ is
locally nilpotent. More precisely, we have the following lemma
(see e.g., \cite[Thm. 16]{Ch} or \cite[Ch. II, Ex. 8]{Ja} for the
case of algebras of finite dimension, \cite{SW} for complete
algebras, and also \cite[2.1]{CD}).

\blem\label{0.2}
\begin{enumerate}[\indent\em (a)]
     \item $A$ is a graded algebra, i.e.\
    $A_\alpha A_\beta\subseteq A_{\alpha+\beta}$ for all $\alpha$,
    $\beta\in k$;
    \item $\delta_s$ and $\delta_n$ are again derivations on $A$;
    \item $\delta_s$ is homogeneous and acts on $A_\alpha$
    via multiplication
    by $\alpha$, so is automatically locally bounded,
    whereas $\delta_n$ is homogeneous of degree $0$ and locally nilpotent.
\end{enumerate}
\elem

\bproof (a) For homogeneous elements $x\in A_\alpha$ and $y\in
A_\beta$ we have
$$
(\delta-(\alpha+\beta)\id)(xy)
=x(\delta-\beta\id)(y)+(\delta-\alpha\id)(x) y\,,
$$
and hence by induction
$$
(\delta-(\alpha+\beta)\id)^n(xy)=\sum_{i=0}^n {n\choose i}
(\delta-\alpha\id)^i(x)\cdot (\delta-\beta\id)^{n-i}(y)\,.
$$
It follows that $(\delta-(\alpha+\beta)\id)^n(xy)$ vanishes for
$n\gg 0$, and so $xy\in A_{\alpha+\beta}$, proving (a).

By definition $\delta_s$ acts via multiplication with  $\alpha$
on $A_\alpha$. Since $A$ is a graded algebra this shows that
$\delta_s$ is a degree $0$ homogeneous derivation on $A$ and so
is $\delta_n=\delta-\delta_s$. Now (b) and (c) follow. \eproof

In the next lemma we consider the set
$$M:=\{\alpha\in k|A_\alpha\ne 0\}\,,$$ and we let $\N M$
and $\Z M$ be the additive subsemigroup and the subgroup of $k$,
respectively, generated by $M$.

\blem\label{0.3} If $A$ is a finitely generated $k$-algebra then
\begin{enumerate}[\indent\em (a)]
     \item also $\N M$ and $\Z M$ are finitely generated, and
     \item $A$ admits an effective action of a torus $(k^*)^r$,
     where $r:=\rk_\Z\Z M$.
\end{enumerate}
\elem

\bproof The proof of (a) is elementary and we omit it. To show
(b) we note that $A=\bigoplus_{\alpha\in M}A_\alpha $ is graded by
the semigroup $\N M$, which is a subsemigroup of $\Z M\cong\Z^r$.
An effective action of $(k^*)^r$ on $A$ is then given by
$$
(\lambda_1,\ldots, \lambda_r).x_\alpha:=\lambda_1^{\alpha_1}\cdots
\lambda_r^{\alpha_r}x_\alpha,\quad x_\alpha\in A_\alpha,
$$
where $(\lambda_1,\ldots, \lambda_r)\in (k^*)^r$ and $\alpha\in
M\subseteq \Z^r$ has components $\alpha_1,\ldots, \alpha_r\in\Z$.
\eproof

\blem\label{0.4} If $\p: A\to A$ is a locally nilpotent
derivation then
\be\label{e1}
\exp(-\p)\delta\exp(\p)=\delta+[\delta, \p]
\ee
for any derivation
$\delta:A\to A$ satisfying
\be\label{e2}
[[\delta,\p],\p]=0\,.
\ee

\elem

\bproof We denote $\Delta_i:=[\delta, \p^i]$. Then (\ref{e2})
says that $[\Delta_1,\p]=0$. By the Jacobi identity we get
\be\label{e3} [\Delta_i, \p]=[\Delta_1,\p^i]\,,\ee  and so by
(\ref{e2}) $[\Delta_i, \p]=0\quad\forall i\ge 1$.

Assuming by induction that, for a given $i$,
\be\label{e4}\Delta_i=i\p^{i-1}\Delta_1\,,\ee we obtain from the
above equalities:
$$\Delta_{i+1}=[\delta,\p^{i+1}]=\Delta_{i}\p+\p^{i}\Delta_1
=i\p^{i-1}\Delta_1\p+\p^i\Delta_1=(i+1)\p^i\Delta_1\,,$$ which
proves (\ref{e4}) for all $i\ge 1$. Thus we get
$$\delta\p^i=\Delta_i+\p^i\delta=\p^i\delta+i\p^{i-1}\Delta_1\,,$$
and consequently
\renewcommand{\arraystretch}{1.6}
$$\begin{array}{rcl}
       \delta\exp(\p)&=&\sum_{i=0}^\infty \delta \p^i/i!\\
        &= & \left(\sum_{i=0}^\infty \p^i/i!\right) \delta+
        \left(\sum_{i=1}^\infty (\p^{i-1}/(i-1)!)\right)\Delta_1\\
         &=&
        \exp(\p) (\delta+[\delta, \p]).
\end{array}
$$
Multiplying with $\exp(-\p)$ from the left gives (\ref{e1}).
\eproof

\section{Main theorem}

In this section we formulate and prove our
main results. For any $\Z$-graded finitely generated $\C$-algebra
$$A=\bigoplus_{i\in\Z}A_i\,,$$
every derivation $\p$ of $A$ has a unique decomposition
$\p=\sum_{i=k}^l\p_i\,,$ where $\p_i: A\to A$ is a
homogeneous derivation of degree $i$. The proof of the following
simple lemma is left to the reader.

\blem\label{0.5} If $\p=\sum_{i=k}^l\p_i$ is locally bounded
    then $\p_l$ and $\p_k$ are also locally bounded.
    In particular, if $l>0$ ($k<0$)
then $\p_l$ (respectively, $\p_k$) is locally nilpotent. \elem

\bsit\label{s7} In the sequel we let $A$ be the coordinate ring of
a normal affine surface with a $\C^*$-action so that
$A=\bigoplus_{i\in\Z}A_i$ is graded. The infinitesimal generator
of this $\C^*$-action is a semisimple derivation $\delta$ on $A$
that acts via $\delta(a)=\deg(a)\cdot a$ for a  homogeneous
element $a\in A$. As was shown in \cite[Proposition 2.4]{FZ2},
for a homogeneous locally nilpotent derivation $\p\neq 0$ the
derivation
$$
\exp(-\p)\delta\exp(\p)
$$
is again semisimple and defines a $\C^*$-action which is, in
general, different from the given one. Conversely, we have the
following result.\esit

\bthm \label{0.6} If $\ML(A)\neq \C$, $\Spec (A)\not\cong
\C^*\times\C^*$ and $\not\cong \A_\C^1\times\C^* $ then every
semisimple derivation $\tdelta$ on $A$ is of the form
$$\tdelta=c\cdot \exp(-\p)\delta\exp(\p),\qquad c\in\C,$$ for some
locally nilpotent derivation $\p$ on $A$. Consequently, any two
effective $\C^*$-actions of $A$, after possibly switching one of
them by the automorphism $\lambda\longmapsto \lambda^{-1}$ of
$\C^*$, are conjugate via an automorphism of $A$ provided by a
$\C_+$-action on $A$ and, moreover, coincide whenever $ML(A)=A$.
\ethm

The latter assertion leads to the following corollary.

\bcor\label{0.7} If a normal affine surface $V=\Spec (A)$,
non-isomorphic to $\C^*\times\C^*$, admits two effective
$\C^*$-actions with infinitesimal generators $\delta, \tilde
\delta$, where $\delta\neq \pm \tilde \delta$, then it also
admits a non-trivial $\C_+$-action.\ecor

To prove Theorem \ref{0.6} we need a few preparations.

\bsit\label{0.8} We suppose below that $\ML(A)\ne \C$, $\Spec
(A)\not\cong \C^*\times\C^*$ and $\not\cong \A_\C^1\times\C^* $.
If $A$ admits a homogeneous locally nilpotent derivation $\p$ of
degree zero then by Lemma 3.8 and Corollary 3.28 in \cite{FZ2}
either $A\cong \C[t,u]$ or $A\cong \C[t,u,u^{-1}]$ with $t\in
A_0,u\in A_d$ homogeneous and $\p=\p/\p t$, which is excluded by
our assumptions. According to Corollary 3.27(i) and Theorem 4.5
from \cite{FZ2}, either all nontrivial homogeneous locally
nilpotent derivations on $A$ are of positive degree, or all of
them are of negative degree. By switching the grading to the
opposite one, if necessary, we may suppose in the sequel that $A$
does not admit a homogeneous locally nilpotent derivation of
degree $\le 0$. \esit

\blem\label{0.9} With the assumptions of \ref{s7} and \ref{0.8},
for every nonzero locally nilpotent derivation $\p$ of $A$ the
following hold.
\begin{enumerate}[\indent\em (a)]
\item $\p$ is a linear combination of commuting homogeneous
locally nilpotent derivations of
strictly positive degrees.
     \item The derivations $[\delta,\p]$ and $\p$ commute.
     \item $\exp(-\p)\delta\exp(\p)=\delta+[\delta,\p]$.
     \item There exists a locally nilpotent derivation $\p'$ on $A$ such
     that $\p=[\delta,\p']$.
\end{enumerate}
\elem

\bproof (a) Let us write $\p$ as a sum of homogeneous derivations
$$
\p=\sum_{i=k}^l\p_i \quad \mbox{with} \quad \p_k, \p_l\ne 0.
$$
Since clearly $\p_k$ and $\p_l$ are again locally nilpotent (see
e.g., \cite{Re}), by our convention in \ref{0.8} above we have
$k\ge l> 0$. Moreover, since $\ML(A)\neq \C$, $\p$ and $\p_k$ are
equivalent, so define equivalent $\A_\C^1$-fibrations, and
$\p=a\p_l$ for some $a\in \Frac (\ker \p_l)$ (see e.g.,
\cite[Lemma 4.5]{FZ2}). It follows that the $\p_i$ are commuting
locally nilpotent derivations, proving (a).

Now (b) follows from (a) since $[\delta,\p]=\sum_{i=k}^li\p_i$,
and (c) follows from Lemma \ref{0.4} by virtue of (b). Finally,
(d) can be deduced by taking $\p':=\sum_{i=k}^l\p_i/i.$ \eproof

\bproof[Proof of Theorem \ref{0.6}] For a semisimple derivation
$\tdelta$ of $A$, we consider its decomposition
$\tdelta=\sum_{i=k}^l \p_i$ into homogeneous components with
$\p_k$, $\p_l\neq 0$. By Lemma \ref{0.5}, if $k<0$ then $\p_k$ is
locally nilpotent, which is excluded by our convention in
\ref{0.8}. Thus $l\ge k\ge 0$.

Now the proof proceeds by induction on $l$. If $l=0$ then
$\tdelta=\p_0$ is semisimple, homogeneous of degree 0 and
commutes with $\delta$. Thus $\delta$ and $\p_0$ are equal up to
a constant factor. Indeed, otherwise $V:=\Spec (A)$ would be a
toric surface non-isomorphic to $\C^*\times\C^*$ or
$\A_\C^1\times\C^* $, hence $\ML(A)=\C$ (see Example 2.8 in
\cite{FZ2}), which contradicts our assumption. Clearly, the
constant factor above is equal to $\pm 1$ as soon as both
$\delta$ and $\p_0$ generate effective $\C^*$-actions.

Assume now that $l>0$. By Lemma \ref{0.9}(d) we find a locally
nilpotent derivation $\p'$ with $[\tdelta, \p']=\p_l$, and so by
Lemma \ref{0.9}(c)
$$
\exp(\p')\tdelta\exp(-\p')=\tdelta-[\tdelta,\p']=\tdelta-\p_l.
$$
Thus $\tdelta-\p_l$ is again semisimple and its homogeneous
components have degrees $\le l-1$. Applying the induction
hypothesis, the result follows. \eproof

\section{Toric surfaces and uniqueness of a DPD-presentation}
\bsit\label{s8} Let $A$ be a normal 2-dimensional $\C$-algebra
with a grading $A=\bigoplus_{i\in\Z}A_i$ associated to an
effective $\C^*$-action. We recall that such a grading admits a
DPD-presentation as follows (see \cite{FZ1}).

{\em Elliptic case:} Here $A_0=\C$, and up to
switching the grading we have $A_-:=\bigoplus_{i<0}A_i=0$. The curve
$C:=\Proj A$ is normal and carries a $\Q$-divisor $D$ of positive
degree unique up to linear equivalence such that

\be\label{eqDPD} A=A_0[D]:= \bigoplus_{i\ge 0} H^0(C,\cO_C(\lfloor
iD \rfloor))u^i\subseteq \Frac(C)[u]\,. \ee \indent{\em Parabolic
case:} Here $A_-=0$, but $A_0$ is 1-dimensional and so defines a
smooth curve $\Spec A_0=\Proj A$. As before $C$ carries a
$\Q$-divisor $D$, now of arbitrary degree and again unique up to
linear equivalence such that (\ref{eqDPD}) holds.

{\em Hyperbolic case:} This case is characterized by $A_+,A_-\ne
0$. The subrings $A_{\ge 0}:=\bigoplus_{i\ge 0}A_i$ and $A_{\le
0}$ are parabolic and as before admit presentations $A_{\ge
0}=A_0[D_+]\subseteq \Frac(C)[u]$ and $A_{\le 0}=A_0[D_-]\subseteq
\Frac(C)[u^{-1}]$, respectively. Thus $A$ is the subring
$$
A=A_0[D_+,D_-]:=A_0[D_+]+A_0[D_-]\subseteq
\Frac(C)[u,u^{-1}].
$$
Moreover by Theorem 4.3 in \cite{FZ1} $D_++D_-\le 0$, and the
pair $(D_+,D_-)$ is determined uniquely by the graded algebra $A$
up to a linear equivalence
$$
(D_+,D_-)\sim (D'_+,D'_-) :\Leftrightarrow (D_+,D_-)=(D'_++\div
f,D'_--\div f)\quad\mbox{with}\quad f\in \Frac (C)^\times\,.
$$
The question arises whether a DPD-presentation is determined
uniquely, up to the linear equivalence as above, an automorphism
of $C$ and in the hyperbolic case by an interchange of $D_+$ and
$D_-$, by the geometry of the surface $V$ alone, disregarding the
choice of a $\C^*$-action. In Corollary \ref{4.1} below we show
that, indeed, this is the case at least for surfaces with a
non-trivial Makar-Limanov invariant.\esit

\bsit\label{s9} Let us recall \cite[Theorem 4.5]{FZ2} that the
Makar-Limanov invariant of a surface $V=\Spec (A)$, where
$A=A_0[D_+,D_-]$ with $D_++D_-\le 0$ and $ D_++D_-\neq 0$, is
trivial if and only if $A_0\cong \C[t]$ and the fractional parts
$\lbrace D_\pm\rbrace$ of $D_\pm$ are zero or are concentrated at
one point, say, $p_\pm\in\A_\C^1$. If still $A_0\cong \C[t]$,
but the second condition holds for precisely one of the divisors
$D_\pm$ then $ML (A)=\C[x]$ for a nonzero homogeneous element
$x\in A$. Otherwise $ML (A)=A$ that is, $A$ does not admit
non-zero locally nilpotent derivations.\esit

  From Theorem \ref{0.6} above and Theorem 4.3 in \cite{FZ1} we
derive the following corollary.

\bcor\label{4.1} For a non-toric normal affine surface
$V=\Spec(A)$ the following hold.
\begin{enumerate}\item[(a)]
If $V$ admits an effective elliptic $\C^*$-action then this
$\C^*$-action is unique up to the automorphism
$\lambda\longmapsto\lambda^{-1}$ of $\C^*$. In particular a
DPD-presentation $A=A_0[D]$ is unique up to linear equivalence of
$D$.
\item[(b)]
If $V$ admits a parabolic $\C^*$-action then every $\C^*$-action
on $V$ is parabolic. Moreover, if $A=A_0[D]=A'_0[D']$ are
DPD-presentations corresponding to two $\C^*$-actions then there
is an isomorphism $\varphi: C\to C'$ of the curves defined by
$A_0$ and $A_0'$, respectively, such that $D$ and $\varphi^*
(D')$ are linearly equivalent.
\item[(c)] If $\ML(A)\neq \C$ and
$$
A=A_0[D_+,D_-]\quad\mbox{and}\quad  A=A_0'[D_+',D_-']
$$
are DPD-presentations of $A$ corresponding to two hyperbolic
$\C^*$-actions on $A$ then there is an isomorphism $\varphi:A_0\to
A_0'$ such that the pair of $\Q$-divisors $(D_+,D_-)$ on the curve
$C=\Spec (A_0)$ is linearly equivalent to one of the pairs
$\varphi^*(D_+',D_-')$ or $\varphi^*(D_-',D_+')$. \end{enumerate}
\ecor

\bproof If in case (a) there is a second $\C^*$-action on $A$ not
related to the first one by an automorphism of $\C^*$ then by
Lemma \ref{0.5} there is a locally nilpotent derivation on $A$
and so by Theorem 3.3 in \cite{FZ2} $V$ is toric contrary to our
assumption.

(b) and (c) follow immediately from the fact that any two
$\C^*$-actions on $V$ are conjugate in $\Aut(V)$ by Theorem
\ref{0.6}. \eproof

\bsit\label{s10} We note that toric surfaces $V=\Spec (A)$ admit
many non-conjugated $\C^*$-actions given by non-conjugated
one-parameter subgroups of the 2-torus $\T^2=\C^*\times \C^*$. It
also has many distinct DPD-presentations, up to permuting $D_+$
and $D_-$ and to linear equivalence. Any pair of divisors
$D_\pm=\mp\frac{e_\pm}{d_\pm}[0]$ with $e_\pm, d_\pm\in \Z$,
$d_+>0,\,d_-<0$, $\gcd(e_\pm,d_\pm)=1$ and $(D_++D_-)(0)<0$
defines a toric surface $V=\Spec( A_0[D_+,D_-])$, where
$A_0\cong\C[t]$. Two such toric surfaces $V=\Spec (A_0[D_+,D_-])$
and $V'=\Spec (A_0[D_+',D_-'])$ are isomorphic if and only if the
sublattices $\Z(e_+,d_+)+\Z(e_-,d_-)$ and
$\Z(e_+',d_+')+\Z(e_-',d_-')$ of $\Z^2$ are equivalent up to the
action of the group of integer matrices with determinant $\pm 1$
(see \cite{De, Gub} for a more general result, or also the proof
of Theorem 4.15 in \cite{FZ1}). According to Theorem \ref{4.2}
below this is the only ambiguity in the choice of a
DPD-presentation for $V$.\esit

In the rest of this section we concentrate on affine toric
surfaces. We remind the reader that any such surface $V=\Spec
(A)$ is isomorphic to $\C^{*2}$, $\A^1_C\times \C^*$ or to
$$
V_{d,e}=\Spec A_{d,e}:=\Spec (\C[X,Y]^{\Z_d})=\A^2_\C/\Z_d,
$$
where the cyclic group
$\Z_d=\langle \zeta\rangle$ of $d$th roots of unity acts on the
polynomial ring $\C[X,Y]$ via
$$
\zeta.X=\zeta X\quad \mbox{and}\quad \zeta.Y= \zeta^eY
$$
for some $e$ satisfying $0\le e< d$ and $\gcd(e,d)=1$. The
standard action of the 2-torus $\T':=\C^{*2}$ on $\A_\C^2=\Spec
(\C[X,Y])$ commutes with the action of $\Z_d$ and so descends to
$V_{d,e}$. Dividing out the kernel $K\cong \Z_d$ gives an
effective action of the 2-torus $\T=\T'/K$ on $V$.

We note that for $ee'\equiv 1\mod d$ the surfaces $V_{d,e}$ and
$V_{d,e'}$ are isomorphic via the map induced by the morphism
$$
\A_\C^2\to\A_\C^2\quad \mbox{with}\quad (x,y)\mapsto (y,x).
$$
This morphism is $\varphi$-equivariant, where
$\varphi$ is the map of $\T'=\C^*\times\C^*$ interchanging the
factors. Moreover, $V_{d,e}$ and
$V_{d',e'}$ are isomorphic as normal surfaces if and only if
\be\label{eqtoric}
d=d'\quad\mbox{and either}\quad
e=e'\quad\mbox{or}\quad ee'\equiv 1\mod d,
\ee
see e.g.\
\cite[Example 2.3]{FZ1}.

Every 1-dimensional subgroup of $\T$ isomorphic to $\C^*$
provides a $\C^*$-action on $V$. Conversely we have the following
result.

\bthm\label{4.2} If $V=\Spec (A)$ is an affine toric surface and if
$\C^*$ acts effectively on $V$ then this action is conjugate to
the action of a subtorus of $\T$.
\ethm

\bproof In the case $V=\T$ this is evident. So we can assume in
the sequel that $V\ne \T$. This $\C^*$-action then defines a
grading $A=\bigoplus_{i\in Z} A_i$ of $A$. As $V$ is toric and
not a torus it admits a $\C_+$-action. Our assertion is an
immediate consequence of the following claim.

\smallskip

{\em Claim.} Either $A=\C[z,v,v^{-1}]$ for homogeneous elements
$z,v$ of $A$, or there is a graded isomorphism $A\cong A_{d,e}$,
where the grading on $A_{d,e}$ is provided by a subtorus of $\T$
as above.

\smallskip

In the case of an  elliptic grading this is just Theorem 3.3 in
\cite{FZ2}. If the grading is parabolic and if there is a
horizontal $\C_+$-action, i.e. the general orbits of the
$\C_+$-action map dominantly to $V/\C^*$, then this is a
consequence of Theorem 3.16 in \cite{FZ2}.

In the case that the grading is parabolic and there is only a
vertical but no horizontal $\C_+$-action, we have $A\cong A_0[u]$
with $\deg u=1$ and $A_0=\C[t,t^{-1}]$. In fact, the curve
$C=\Spec (A_0)$ is either $\A^1_C$ or $\C^*$, since the open
orbit of the torus action on $V$ maps dominantly onto $C$. If
$C\cong\A^1_\C$ then $A=A_0[D]$ with $A_0\cong \C[t]$ and the
fractional part of $D$ is supported on at least two points, since
by our assumption there is no horizontal $\C_+$-action on $V$
(cf.\ \cite[Theorem 3.16]{FZ2}). By Proposition 3.8 in \cite{FZ1}
this would imply that $V$ has at least two singular points, which
is impossible.

Therefore $C\cong \C^*$ and so $A_0\cong \C[t,t^{-1}]$. As
there is no dominant morphism $\A^2_\C\to \C^*$ and thus also no
dominant morphism $V_{d,e}\to\C^*$, we have $V\cong
\A^1_\C\times\C^*$. In particular $V$ is smooth, and $A=A_0[D]$ with
$D$ being an integral (so, principal) divisor. Thus $D\sim 0$ proving
the claim in this case.

In the remaining case the grading on $A$ is hyperbolic so that
$A=A_0[D_+,D_-]$. Since $V$ admits a $\C_+$-action, we have
necessarily $A_0\cong \C[t]$ and $C=\A^1_\C$, see \cite[Corollary
3.23]{FZ2}. The Picard group of a toric surface is trivial and so
by Corollary 4.24 in \cite{FZ1} we have $l\le 1$, where
$$
l:=\mbox{card}\ \{p\in C :
D_+(p)+D_-(p)<0\}.
$$
If $l=0$ then $D_+=-D_-$ and so by Remark 4.5 in \cite{FZ1}, $A$
contains a unit of nonzero degree. As $V$ admits a $\C_+$-action
then by Corollary 3.27 in \cite{FZ2} we have $A=\C[z,v,v^{-1}]$
with homogeneous elements $z,v\in A$, which proves our claim in
this case.

If $l=1$ then by Corollary 4.24 in \cite{FZ1} the orbit map $V\to
C=V//\C^*$ has no irreducible multiple fiber and the divisor
$D_++D_-$ is concentrated in one point, say, $p\in\A_\C^1$. Since
$D_\pm(a)\in\Z$ and $D_+(a)+D_-(a)=0$ for every point $a\neq p$,
the pair $(D_+,D_-)$ is equivalent to a pair
$$
\left( -\frac{e_+}{d_+}[p],\frac{e_-}{d_-}[p]\right)\,.
$$
As in the proof of Theorem 4.15 in \cite{FZ1} there exist
integers $d,e$ with $0\le e<d$ and $\gcd(d,e)=1$ so that $A\cong
A_{d,e}$ as graded rings, where the grading on $A_{d,e}$ is
defined by a subgroup of $\T$ isomorphic to $\C^*$. This finishes
the proof. \eproof

As a particular case we obtain the following classical
result.

\bcor \label{4.3}
{\em (Gutwirth \cite{Gut})} Every $\C^*$-action on $\A_\C^2$
is linearizable. \ecor

\brems\label{4.4} 1. Any two effective actions of a 2-torus on a
normal affine surface $V$ are conjugate. This follows from the
fact that any toric surface is equivariantly isomorphic to one of
the surfaces $\C^{*2}$, $\A^1_\C\times\C^*$ or $V_{d,e}$, and two
such surfaces are isomorphic as abstract surfaces if and only if
(\ref{eqtoric}) holds.

2. Every normal affine surface $V=\Spec (A)$ with an elliptic or
parabolic $\C^*$-action and with a trivial Makar-Limanov
invariant is toric (see Theorems 3.3 and 3.14 in \cite{FZ2}).
This yields that there can be at most one parabolic
DPD-presentation $A=A_0[D]$ on an affine normal surface, up to the
equivalence described in Corollary \ref{4.1}(b). \erems

\section{Homogeneous affine surfaces}

In this
section we show the following result.

\bthm\label{5.1}
Let $V$ be
one of the surfaces $\PP^1\times \PP^1\backslash\Delta$ or
$\PP^2\backslash Q$, where $\Delta$ is the diagonal in
$\PP^1\times\PP^1$ and $Q$ is a nonsingular quadric in $\PP^2$. Then
any two $\C^*$-actions on $V$ are conjugate in $\Aut(V)$.
\ethm

In case $V=\PP^2\backslash Q$ this is a result of \cite{DG},
see also \cite{Do}. This theorem follows immediately from
the following two more general statements.

\bprop\label{5.2}
Every  smooth Gorenstein
$\C^*$-surface $V=\Spec(A)$ with $\Cl (V)\cong\Z$ and $\ML(V)=\C$ is
isomorphic to $V\cong \PP^1\times \PP^1\backslash\Delta$, where
$\Delta$ is the diagonal. Moreover, there is a graded isomorphism
$A\cong A_0[D_+,D_-]$, where $A_0=\C[t]$, $D_+=0$ and
$D_-=-[1]-[-1]$.
\eprop

\bproof An elliptic or parabolic $\C^*$-surface with trivial
Makar-Limanov invariant is toric by Theorems 3.3 and 3.16 from
\cite{FZ2}, and so its Picard group vanishes. Thus under our
assumptions the $\C^*$-action on $V$ is necessarily hyperbolic
and provides a DPD-presentation $A=A_0[D_+,D_-]$ with $A_0=\C[t]$
and a pair of $\Q$-divisors $D_\pm$ on the affine line
$\A^1_\C=\Spec \C[t]$ such that $D_++D_-\le0$, see Corollary 3.23
in \cite{FZ2}. As the divisor class group of $A$ is isomorphic to
$\Z$ by Theorem 4.22 in \cite{FZ1} we have
$$
l=\mbox{card}\ \{p\in\A^1_\C: D_+(p)+D_-(p)<0\}=2.
$$
Thus there are unique points $p_1,p_2\in \A^1_\C$ with
$D_+(p_i)+D_-(p_i)<0$ for $i=1,2$. Moreover there can be no
multiple fiber of the projection $V\to \A^1_\C$. Indeed, by {\em
loc.\ cit.} every such fiber contributes to the torsion of
$\Cl(A)$. Therefore we may assume that
$$
D_+(p_i)=-\frac{e^+_i}{m_i^+}\quad\mbox{and}\quad
D_-(p_i)=\frac{e_i^-}{m_i^-}
$$
with $m_i^+> 0$ and $m_i^-<0$ and $\gcd(e_i^\pm,m_i^\pm)=1$ for
$i=1,2$. Since $V$ is smooth, we have
\be
\label{eq2}
\left|\begin{array}{cc}
     e_i^+  & m_i^+   \\
     e_i^-
&  m_i^-
\end{array}\right|=-1,
\ee see \cite[Theorem 4.15]{FZ1}. In this case using again
Theorem 4.22 from \cite{FZ1} the divisor class group is generated
by the orbit closures $O_i^\pm$ over the points $p_i$ modulo the
relations

$$
M_1=M_2=0\quad\mbox{and} \quad
E_1+E_2=0,
$$
where
$$
M_i:=m_i^+[O_i^+]-m_i^-[O_i^-]\quad
\mbox{and}\quad
E_i:=e_i^+[O_i^+]-e_i^-[O_i^-].
$$
By assumption the canonical divisor of $V$ is trivial, and by
Corollary 4.25 in \cite{FZ1} it is given by
$$
K_V=\sum_{i=1}^2
((m_i^+-1)[O_i^+]+(-m_i^--1)[O_i^-])\sim
-\sum_{i=1}^2
([O_i^+]+[O_i^-]).
$$
Thus the divisor $K=-\sum ([O_i^+]+[O_i^-])$ is contained in the
subgroup generated by $M_1$, $M_2$ and $E_1+E_2$, and so the
determinant $\det(K,M_1,M_2,E_1+E_2)$ vanishes. Using (\ref{eq2})
this leads to the relation \be\label{eq3} m_1^++m_1^-=m_2^++m_2^-.
\ee Since the Makar-Limanov invariant of $V$ is trivial, the
fractional part of each of the divisors $D_\pm$ is concentrated
in one (possibly the same) point, see \ref{s9}. Thus we may
assume that $m_1^+=1$, and $m_i^-=-1$ for at least one value
$i=1,2$. In the case $m_1^-=-1$ equation (\ref{eq3}) shows that
$m_2^++m_2^-=0$. By (\ref{eq2}) this yields $m_2^+=-m_2^-=1$.
Similarly, if $m_2^-=-1$ then the right hand side in (\ref{eq3})
is $\ge 0$ whereas the term on the left is $\le 0$. This forces
$m_2^+=1$ and $m_1^-=-1$.

Hence the divisors $D_\pm$ are necessarily integral. Replacing
them by an equivalent pair of divisors we can achieve that
$D_+=0$ and $D_-=-[p_1]-[p_2]$. After performing an automorphism
of $\A^1_\C$ we can also assume that $p_1=1$ and $p_2=-1$ and so
the DPD-presentation has the required form. Now the isomorphism
of $V$ and $\PP^1\times\PP^1\backslash \Delta$ follows from
Example 5.1 in \cite{FZ1}. \eproof

\bprop\label{5.3}
A smooth
$\C^*$-surface $V=\Spec A$ with $\Pic (V)\cong\Z_2$ and $\ML(V)=\C$
is isomorphic to $\PP^2\backslash Q$, where $Q$ is a smooth quadric
in $\PP^2$. Moreover there is a graded isomorphism $A\cong
A_0[D_+,D_-]$, where $A_0=\C[t]$, $D_+=\frac{1}{2}[0]$ and
$D_-=-\frac{1}{2}[0]-[1]$.
\eprop

\bproof With the same arguments as in the proof of Proposition
\ref{5.2} the given $\C^*$-action is necessarily hyperbolic and
$V$ admits a DPD-presentation $V=\Spec (A_0[D_+,D_-])$ with
$A_0=\C[t]$. As $V$ is smooth, by Theorem 4.22 in \cite{FZ1} the
condition $\Pic(V)=\Z_2$ gives that $l=1$ and there is only one
irreducible multiple fiber of multiplicity 2. Thus we may suppose
that this multiple fiber lies over $0$ and then the fractional
parts are
$$
\{D_\pm(0)\}=1/2.
$$
The fractional parts of $D_\pm$ at all other points must be zero
since $\ML(A)=\C$ (see Theorem 4.5 in \cite{FZ2}). Hence after
passing to an equivalent pair of divisors, we may suppose that
$D_+= \frac{1}{2}[0]$ and then $D_-=-\frac{1}{2}[0]-a[p]$ for some
$a\in\N$ and $p\in \A^1_\C$. After performing an automorphism of
$\A^1_\C$, if nedeed, we may also suppose that $p=1$. As $V$ is
smooth we obtain from Theorem 4.15 in \cite{FZ1} that necessarily
$a=1$. Thus
$$
A=A_0[D_+,D_-]\quad \mbox{with}\quad
D_+=\frac{1}{2}[0]\quad
\mbox{and}\quad
D_-=-\frac{1}{2}[0]-[1].
$$
This proves the second assertion. The first one follows by
comparing with Example 5.1 from \cite{FZ2}. \eproof

\brem\label{conjecture}
If  a
torus of dimension $k$ acts effectively on a variety $V$ of dimension
$n$, then $k\le n$. Thus the rank of a reductive group $G$ acting
effectively on $V$ is bounded by $n$. As every increasing chain of
connected reductive groups with bounded rank becomes stationary, this
shows the following fact: {\em every connected reductive subgroup of
$\Aut(V)$ is contained in a maximal one.}

In the Introduction we
posed the question whether any two maximal connected reductive
subgroups of the automorphism group of a normal affine surface are
conjugate. We note that this would follow if one could prove that

\noindent $(*)$  any two effective $\C^*$-actions on a non-toric
surface $V$ with $\ML(V)=\C$ are conjugate.

In fact, if $V$ admits
no  $\SL_2$-action then every maximal connected reductive subgroup of
$\Aut(V)$ is a torus. In this case the result follows from $(*)$,
Theorems \ref{0.6}, \ref{4.2} and Remark \ref{4.4}(1).

If $V$ admits
an action of $\SL_2$ then by the Theorem of Gizatullin and Popov
mentioned in the Introduction $V$ is one of the surfaces $\PP^1\times
\PP^1\backslash \Delta$,  $\PP^2\backslash Q$, or $V_{d,1}$ ($d\ge
1$).

If $V$ is one of the surfaces $\PP^1\times \PP^1\backslash
\Delta$ or $\PP^2\backslash Q$ then the standard actions of
$\PGL_2$ on $V$ cannot be extended to a larger connected
reductive group, since otherwise there would be an action of a
2-torus on $V$. Thus in these cases the maximal connected
reductive subgroups of $\Aut(V)$ are isomorphic to $\PGL_2$ and,
moreover, any two such subgroups are conjugate in $\Aut(V)$ by
Proposition 4.14 in \cite{FZ2}.

Similarly, if $V=V_{d,1}$ then $\GL_2/\Z_d$ is a maximal
connected reductive subgroup of $\Aut(V)$. Given another maximal
connected reductive subgroup $G$ of $\Aut(V)$ we may suppose by
Remark \ref{4.4}(1) that its maximal torus is equal to the standard
one in $\GL_2/\Z_d$. Now it is easy to see that $G$ and $\GL_2/\Z_d$
are equal (alternatively one can apply Lemma 4.17 in
\cite{FZ2}).
\erem

 From this remark and Theorems \ref{0.6}, \ref{4.2} and \ref{5.1}
we deduce the following corollary.

\bcor\label{fc} For a normal affine surface $V$, any two maximal
connected reductive subgroups of the automorphism group $\Aut(V)$
are conjugate in $\Aut(V)$ except possibly in the case where $V$
is non-toric, has a trivial Makar-Limanov invariant and is not
isomorphic to one of the surfaces $\PP^1\times \PP^1\backslash
\Delta$ or $\PP^2\backslash Q$. \ecor

The forthcoming paper \cite{FKZ} will be devoted to this
remaining case.

{\it Added in proofs.} We are grateful to Peter Russell who showed
us an example, given any $k\ge 2$, of a non-toric smooth affine
surface with a trivial Makar-Limanov invariant that admits $k$
mutually non-conjugated $\C^*$-subgroups in the automorphism
group. This gives, in general, a negative answer to the question
in the Introduction. Such a surface appears as the complement in a
Hirzebruch surface of a section with selfintersection number
$k+1$. These surfaces were studied in \cite[II]{DG}.

We are also grateful to J\"urgen Hausen who informed us, after
appearing of the first e-print version of this paper, that in
\cite{BeHau} an $n$-dimensional generalization of our Theorem
\ref{4.2} had been proven. Namely, it was shown that an effective
action of the torus $\T^{n-1}$ on an $n$-dimensional affine toric
variety $V$ is conjugate in the automorphism group $\Aut (V)$ to
a subgroup of the standard action of $\T^{n}$ on $V$.


\begin{thebibliography}{KMMH}

\bibitem[BB]{BB}
A.\
Bialynicki-Birula, {\em Remarks on the action of an algebraic torus
on $k^n$. I}, Bull.\ Acad.\ Polon.\ Sci.\ Math.\ Astr., Phys.\ XIV
(1967), 177--188; {\it II, ibid.}, XV (1967), 123--125.

\bibitem[BML]{BML}
T.\ Bandman, L.\ Makar-Limanov, {\em Affine surfaces with
$AK(S)=\C$}, Michigan Math.\ J.\ 49 (2001), 567--582.

\bibitem[BaHab]{BaHab} H.\ Bass, W.\ Haboush, {\em
Linearizing certain reductive group actions}, Trans.\ Amer.\
Math.\ Soc.\ 292 (1985), 463--482.

\bibitem[BeHau]{BeHau} F. Berchtold, J. Hausen, {\em Demushkin's theorem in codimension
one}, Math.\ Z.\ 244 (2003), 697--703.

\bibitem[Be]{Be} J.\ Bertin,
{\em  Pinceaux de droites et automorphismes des surfaces affines},
J.\ Reine Angew.\ Math.\ 341 (1983), 32--53.

\bibitem[Ch]{Ch} C.\ Chevalley, {\em Th\'eorie des groupes de Lie. Tome II.
Groupes alg\'ebriques}, Hermann-Cie.\, Paris, 1951.

\bibitem[CD]{CD} A.\ M.\ Cohen, J.\ Draisma, {\em From Lie algebras
of vector fields
to algebraic group actions}, Transform. Groups 8 (2003), 51--68.

\bibitem[Da]{Da} D.\ Daigle, {\em On locally nilpotent derivations of
$k[X,Y,Z] / (XY-P(Z))$}, J.\ Pure Appl.\ Algebra 181 (2003),
181--208.

\bibitem[DG]{DG}
V.\ I.\ Danilov, M.\ H.\ Gizatullin, {\em Automorphisms of affine
surfaces. I}, Math. USSR Izv. 9 (1975), 493--534; {\it II, ibid.} 11
(1977), 51--98.

\bibitem[De]{De} A.\ S.\ Demushkin, {\em Combinatorial invariance of toric
singularities}, Vestnik Moskov.\ Univ.\ Ser.\ I Mat.\ Mekh.\
1982, 80--87.

\bibitem[Do]{Do} M.\ Doebeli, {\em Linear models for reductive group
actions on affine quadrics}, Bull.\ Soc.\ Math.\ France 122
(1994), 505--531.

\bibitem[Du]{Du} A.\ Dubouloz, {\em Completions of
normal affine surfaces with a trivial Makar-Limanov invariant},
Pr\'epublication de l'Institut Fourier n$^\circ $ 579 (2002).


\bibitem[FKZ]{FKZ} H.\ Flenner, S.\ Kaliman, M.\ Zaidenberg, {\em
On the uniqueness of $\C^*$-actions on normal affine surfaces}.
{\bf II} ({\it in preparation}).

\bibitem[FZ$_1$]{FZ1}
H.\ Flenner, M.\ Zaidenberg, {\em Normal affine surfaces with $\bf
C^*$-actions}, Osaka J.\ Math.\ 40:4, 2003,  981--1009.

\bibitem[FZ$_2$]{FZ2}
H.\ Flenner, M.\ Zaidenberg, {\em Locally nilpotent derivations on
affine surfaces with a ${\bf C}^*$-action.} Pr\'epublication de
l'Institut Fourier de Math\'ematiques, 638, Grenoble 2004,  35p.
math.AG/0403215.

\bibitem[Gi]{Gi} M.\ H.\ Gizatullin,
{\em I. Affine surfaces that are quasihomogeneous with respect to
an algebraic group}, Math. USSR Izv. 5 (1971), 754--769; {\em II.
Quasihomogeneous affine surfaces}, {\it ibid.} 1057--1081.

\bibitem[Gub]{Gub} J.\ Dzh.\ Gubeladze, {\em I. The isomorphism problem
for monoid rings of
rank $2$ monoids}, J.\ Pure Appl.\  Algebra 105 (1995), 17--51;
{\em II. The isomorphism problem for commutative monoid rings},
{\it ibid.} 129 (1998), 35--65.

\bibitem[Gut]{Gut} A.\ Gutwirth,
{\em The action of an algebraic torus on the affine plane}, Trans.
Amer.\ Math.\ Soc.\ 105 (1962), 407--414.

\bibitem[HM]{HM} H.\ Hauser, G.\ M\"uller, {\em Algebraic singularities
have maximal
reductive automorphism groups}, Nagoya Math.\ J.\ 113 (1989),
181--186.

\bibitem[Ja]{Ja} N.\ Jacobson, {\em Lie algebras}, Dover Publ.\ Inc.\, NY 1979.

\bibitem[KKMR]{KKMLR} S.\ Kaliman, M.\ Koras, L.\ Makar-Limanov, P.\ Russell,
{\em $\mathbf C\sp *$-actions on $\mathbf C\sp 3$ are linearizable},
Electron.\ Res.\ Announc.\ Amer.\ Math.\ Soc.\ 3 (1997), 63--71.

\bibitem[Kam]{Kam} T.\ Kambayashi, {\em Automorphism group of a
polynomial ring and
algebraic group action on an affine space}, J.\ Algebra 60
(1979), 439--451.

\bibitem[Kn]{Kn} F.\ Knop, {\em
Nichtlinearisierbare Operationen halbeinfacher Gruppen auf affinen
Räumen}, Invent.\ Math.\ 105 (1991), 217--220.

\bibitem[KP]{KP} H.\ Kraft, V.\ L.\ Popov, {\em Semisimple group actions on
the three-dimensional affine space are linear}, Comment.\ Math.\
Helv.\ 60 (1985), 466--479.

\bibitem[KS]{KS} H.\ Kraft, G.\ W.\ Schwarz, {\em Reductive
group actions with one-dimensional quotient}, Inst.\ Hautes
\'Etudes Sci.\ Publ.\ Math.\ 76 (1992), 1--97.

\bibitem[ML$_2$]{ML2} L.\ Makar-Limanov,
{\em  On the group of automorphisms of a surface $x^ny=P(z)$},
Israel J.\ Math.\ 121 (2001), 113--123.

\bibitem[MMP]{MMP} M.\
Masuda, L.\ Moser-Jauslin, T.\ Petrie, {\em Equivariant algebraic
vector bundles over representations of reductive groups:
applications}, Proc.\ Nat.\ Acad.\ Sci.\ U.S.A.\ 88 (1991), no.\ 20,
9065--9066.

\bibitem[Pa]{Pa} D.\ I.\ Panyushev, {\em Semisimple
groups of automorphisms of a four-dimensional affine space}, Izv.\
Akad.\ Nauk SSSR Ser.\ Mat.\ 47 (1983), 881--894. A correction:
{\it ibid.} 49 (1985), 1336--1337.

\bibitem[Po]{Po}
V.\ L.\ Popov,
{\em On polynomial automorphisms of affine spaces}, Izv.\ Math.\ 65
(2001), 569--587.

\bibitem[Re]{Re} R.\ Rentschler,
{\em Op\'erations du groupe additif sur le plane affine}, C.\ R.\
Acad.\ Sci.\ {\bf 267} (1968), 384-387.

\bibitem[SW]{SW} G.\ Scheja, H.\ Wiebe, {\em
\"Uber Derivationen von lokalen analytischen Algebren}, Symposia
Mathematica XI. Academic Press, London, 1973, 161--192.

\bibitem[Sc]{Sc} G.\ W.\ Schwarz, {\em Exotic algebraic group
actions}, C.\ R.\ Acad.\ Sci.\ Paris S\'er. I Math. 309 (1989),
89--94.

\end{thebibliography}
\end{document}